\documentclass[a4j,11pt]{article}

\usepackage{amsthm}
\usepackage{amsmath}
\usepackage{amssymb}
\usepackage{amsfonts}
\usepackage{wrapfig}
\usepackage[pdftex]{graphicx,color}
\usepackage{comment}

\theoremstyle{plain}

\newtheorem{thm}{Theorem}
\newtheorem{lem}{Lemma}

\newcommand{\argmax}{\mathop{\rm arg~max}\limits}

\title{A New Model Variance Estimator for an Area Level Small Area Model to Solve Multiple Problems Simultaneously}
\author{Masayo Yoshimori Hirose \\{The Institute of Statistical Mathematics}\\
    and \\
    Partha Lahiri\\
    Joint Program in Survey Methodology, University of Maryland, \\College Park, U.S.A.
 }
\date{\hfill}

\begin{document}

\maketitle

\begin{abstract}
 The two-level normal hierarchical model (NHM) has played a critical role in the theory of {\it small area estimation} (SAE), one of the growing areas in statistics with numerous applications in different disciplines. In this paper, we address major well-known shortcomings associated with the empirical best linear unbiased prediction  (EBLUP) of a small area mean and its mean squared error  (MSE) estimation by considering an appropriate  model variance estimator that satisfies multiple properties. The proposed model variance estimator simultaneously (i) improves on the estimation of the related shrinkage factors, (ii) protects EBLUP from the common overshrinkage problem, (iii) avoids complex bias correction in generating strictly positive second-order unbiased mean square error (MSE) estimator either by the Taylor series or single parametric bootstrap method.  The idea of achieving multiple desirable properties in an EBLUP method through a suitably devised model variance estimator is the first of its kind  and holds promise in providing good inferences for small area means under the classical linear mixed model prediction framework.  The proposed methodology is also evaluated using a Monte Carlo simulation study and real data analysis.
\end{abstract}


\section{Introduction}
\label{sec:intro}



Planning and evaluation of government programs usually requires access to a wide range of national and sub-national socio-economic, environment and health related statistics. There is, however, a growing need for statistics relating to much smaller geographical areas where data are too sparse to support the sort of standard estimation methods typically employed at the national level. These small area official statistics are routinely used for a variety of purposes, including assessing economic well-being of a nation, making public policies, and allocating funds in various government programs. In this context, the term small area typically refers to a sub-population for which reliable statistics of interest cannot be produced using the limited area specific data available from the primary data source.

With the availability alternative data sources such as survey data, administrative and census records, different governmental agencies are now exploring ways to combine information from different data sources in order to produce reliable small area statistics.  A common practice is to use a statistical model, usually a mixed model, and an efficient statistical methodology such as Bayesian or EBLUP for combining information from multiple databases.  Such a strategy generally improves on estimation for a domain with small or no sample from the primary data source.    We refer to the book by Rao and Molina (2015) for a comprehensive recent account of small area estimation literature.

Both classical and Bayesian methods and theories have been developed using the following widely applied two-level Normal hierarchical model:
\vskip .2in
\noindent{\bf A Two-Level Normal Hierarchical Model (NHM)}
\begin{description}
  \item {Level 1 (sampling model):}  $y_i|\theta _i
\stackrel{\mathrm{ind}}{\sim}
  N(\theta_i,D_i)$;
  \item {Level 2 (linking model):} $\theta_i
\stackrel{\mathrm{ind}}{\sim}
  N(x_i^{\prime}\beta, A),$
\end{description}
for $i=1,\cdots,m.$
\vskip .2in
In the above model, level 1 is used to account for the sampling distribution of unbiased estimates $y_i$.  For example, $y_i$ could be a sample mean based on  $n_i$ observations taken from the $i$th population (e.g., a small geographic area, a hospital or a school.)  As in other papers on the  NHM (e.g., Efron and Morris 1973, 75; Fay and Herriot 1979; Morris 1983; Datta, Rao, and Smith, 2005), we assume that the sampling variances $D_i$ are known, in order to concentrate on the main issues.  The assumption of known sampling variances $D_i$ often follows from the asymptotic variances of transformed direct estimates (Efron and Morris 1975; Carter and Rolph 1974) and/or from empirical variance modeling (Fay and Herriot 1979, Bell and Otto 1995).

Level 2 links the random effects $\theta_i$ to a vector of $p$ known auxiliary variables $x_i=(x_{i1},\cdots,x_{ip})^{\prime}$, often obtained from various alternative data sources (e.g., administrative records, severity index for a hospital, school register, etc.). The parameters $\beta$ and $A$ of the linking model, commonly referred to as hyperparameters, are generally unknown and are estimated from the available data. We assume that $\beta\in R^p,$ the $p$-dimensional Euclidian space, and $A\in [0,\infty)$.

The NHM model can be viewed as the following simple linear mixed model:
\begin{eqnarray}
y_i=\theta_i+e_i=x_i^{\prime}\beta+v_i+e_i, \ (i=1,\ldots,m),\label{FH}
\end{eqnarray}
where $\{v_1\ldots,v_m\}$ and $\{e_1,\ldots, e_m\}$ are independent with $v_i{\sim}N(0,A)$ and $e_i{\sim}N(0,D_i)$; $x_i$ is a $p$-dimensional vector of known auxiliary variables; $\beta\in R^p$ is a $p$-dimensional vector of unknown regression coefficients; $A\in [0,\infty)$ is an unknown variance component; $D_i>0$ is the known sampling variance of $y_i\;(i=1,\cdots,m)$.

NHM is particularly effective in combining different sources of information and explaining different sources
of errors. Some earlier applications of NHM include the estimation of: (i) false alarm probabilities in
New York city (Carter and Rolph 1974), (ii) the batting averages of major league baseball players (Efron and Morris 1975), and (iii) prevalence of toxoplasmosis in El Salvador (Efron and Morris 1975).

Since the publication of the landmark paper by Fay and Herriot (with 971 google citation to date), the NHM, commonly known as the {\it Fay-Herriot (FH) model} in the small area research community, has been  extensively used in developing small area estimation theory and in a wide range of applications.  In a small area estimation setting, NHM or the FH  was used: to estimate poverty rates for the US states, counties, and
school districts  (Citro and Kalton 2000) and Chilean municipalities (Casas-Codero et al. 2015), and to estimate proportions at the lowest level of literacy for states and counties (Mohadjer et al. 2007).

%

The MSE of a given predictor $\hat{\theta}_i$ of  $\theta_i$ is defined as $M_i(\hat{\theta}_i)=E(\hat{\theta}_i-\theta_i)^2$, where the expectation is with respect to the joint distribution of $y=(y_1,\cdots,y_m)^{\prime}$ and $\theta=(\theta_1,\cdots,\theta_m)^{\prime}$ under the Fay--Herriot model (\ref{FH}).
The best linear unbiased predictor (BLUP)  $\hat{\theta}_i^{BLUP}$ of $\theta_i$, which minimizes  $M_i(\hat{\theta}_i)$ among all linear unbiased predictors $\hat\theta_i$, is given by:
$$\hat{\theta}_i^{BLUP} (A)=(1-B_i)y_i+B_i x^{\prime}_i\hat{\beta}(A),$$
where $B_i\equiv B_i(A)=D_i/(A+D_i)$ is the shrinkage factor and  $\hat{\beta}(A)=(X^{\prime}{V}^{-1}X)^{-1}X^{\prime}{V}^{-1}y$ is the weighted least square estimator of $\beta$ when $A$ is known.  Here we use the following notation: $X^{\prime}=(x_1,\cdots,x_m),$ a $p\times m$  matrix of known auxiliary variables; $V=\mbox{diag}(A+D_1,\cdots,A+D_m),$ a $m\times m$ diagonal matrix. By plugging in an estimator $\hat A$ for $A$ (e.g., ML, REML, ANOVA) in the BLUP, one gets an empirical BLUP (EBLUP): $\hat{\theta}_i^{EB}\equiv \hat {\theta}_i^{BLUP} (\hat A)$.


 In the context of an empirical Bayesian approach, Morris (1983) noted that for making inferences about $\theta_i$, estimation of $B_i$ is more important than that of $A$ because the posterior means and variances of $\theta_i$ are linear in $B_i$, not in $A$.  He also noted that, even if an exact unbiased estimator of $A$ is plugged in $B_i\equiv B_i(A)$, one may estimate $B_i$ with large bias.  For that reason, to motivate the James-Stein estimator of $\theta_i$, Efron and Morris (1973) used an exact unbiased estimator of $B$ and {\it not} maximum likelihood estimator of $A$.  For small $m$, maximum likelihood estimator of $A$ (even with the REML correction) frequently produces estimate of $A$ at the boundary (that is, 0) resulting in $B_i=1$ for all $i$,  even when some of the true $B_i$ are not close to 1.  This causes an overshrinkage problem in EBLUP.  That is, for each $i$, EBLUP of $\theta_i$ reduces to the regression estimator. To overcome the overshrinkage problem, Morris (1983) suggested the fraction $(m-p-2)/(m-1)$ when estimator of $B_i$ is 1.  Li and Lahiri (2010) and Yoshimori and Lahiri (2014) avoided the overshrinkage problem by considering strictly positive consistent estimators of $A$, but did not devise their estimators of $A$ to obtain nearly accurate estimator of $B_i$; that is, biases of their estimators of $B_i$, like all other existing estimators (e.g., ML or REML), are of the order $O(m^{-1})$ and not $o(m^{-1})$.  This is an important research gap, which we will fill in this paper.

An estimator $\hat M_i(\hat\theta_i^{EB})$ of $M_i(\hat\theta_i^{EB})$ is called second-order unbiased if $E[\hat M_i(\hat\theta_i^{EB})]  =M_i(\hat{\theta}_i^{EB})+o(m^{-1}),$  for large $m$, under suitable regularity conditions.  Let $M_{i;approx}(A)$ be a second-order approximation to $M_i(\hat{\theta}_i^{EB}).$  That is, $M_i(\hat{\theta}_i^{EB})=M_{i;approx}(A)+o(m^{-1}),$ for large $m$, under regularity conditions.  Prasad and Rao (1990) proposed a second-order unbiased estimator of $M_i(\hat\theta_{i;MOM}^{EB})$, where $\hat{\theta}_{i;MOM}^{EB}$ is EBLUP of $\theta_i$ when method-of-moments (MOM) estimator $\hat A_{MOM}$ of $A$ is used.  They noticed that the simple plugged-in estimator  $M_{i;approx}(\hat A_{MOM})$ is not second-order unbiased estimator of $M_i(\hat{\theta}_{i;MOM}^{EB})$.  They showed that
 $$E[M_{i;approx}(\hat A_{MOM})]=M_i[\hat{\theta}_i^{EB}(\hat A_{MOM})]+O(m^{-1}),$$
 for large $m$, under regularity conditions. In fact, $M_{i;approx}(\hat A)$ is not second-order unbiased estimator of $M_i(\hat{\theta}_i^{EB})$ for any variance component estimators proposed in the literature.
 Bias correction is usually applied to achieve second-order unbiasedness. However, some bias-correction can even yield negative estimates of MSE.   See Jiang (2007) and Molina and Rao (2015) for further discussions.


\indent Mimicking a Bayesian hyperprior calculation, Laird and Louis (1987) introduced a parametric bootstrap method for measuring uncertainty of an empirical Bayes estimator.  While their point estimator is identical to EBLUP, their measure of uncertainty has more of a Bayesian flavor rather than MSE.  Butar (1997) [see also Butar and Lahiri 2003] was the first to introduce parametric bootstrap method to produce a second-order unbiased MSE estimator in the small area estimation context.  Since Butar's work, a number of papers on parametric bootstrap MSE estimation methods appeared in the SAE literature; see Pfeffermann and Glickman (2004), Chatterjee and Lahiri (2007); Hall and Maiti (2006); Pfefferman and Correra (2012).  Some of them are the second-order unbiased but not strictly positive.  Some adjustments were proposed to make the second-order unbiased double parametric bootstrap MSE estimators strictly positive, but adjusted MSE estimators were not claimed to have the dual property of second-order unbiasedness and strict positivity.  As pointed out in Jiang et al. (2016), a proof is not at all trivial and it is not even clear if the adjustments for positivity retain the second-order unbiasedness of the MSE estimators.

In this paper, we focus on the estimation of two important area-specific functions of $A$ --- the shrinkage factor $B_i$ and the MSE of the EBLUP $M_i(\hat\theta_i^{EB})$.
We propose a single area specific estimator of $A$, say $\hat A_i,$ that simultaneously satisfies the following multiple desirable properties under certain mild regularity conditions:
\begin{description}
\item Property 1: Obtain a second-order unbiased estimator of $B_i$, that is, $E(\hat B_i)=B_i+o(m^{-1})$, among the class of estimators of $B_i$
with identical variance, up to the order $O(m^{-1})$, where $\hat B_i=D_i/(\hat A_i+D_i)$.
\item Property 2: $0<\mbox{inf}_{m\ge 1}\hat B_i\le \mbox{sup}_{m\ge 1}\hat B_i<1$.  That is, it protects  EBLUP from overshrinking to the regression estimator, a common problem encountered in the EB method;
\item Property 3:  Obtain  second-order unbiased Taylor series MSE estimator of EBLUP without any bias correction; that is, $E[M_{i;approx}(\hat A_i)]=M_i(\hat{\theta}_i^{EB})+o(m^{-1}).$
\item Property 4:  Produce a strictly positive second-order unbiased single parametric bootstrap MSE estimator without any bias-correction.
\end{description}
Note that the variance component in the FH model (\ref{FH}) is not area specific, but to satisfy the above properties simultaneously for a given area, we propose an area specific estimator of $A$.  This introduces an area specific bias, but interestingly the order of bias is $O(m^{-1})$,  same as the bias of the ML estimator of $A$ but higher than that of REML in the higher-order asymptotic sense.  This seems to be a reasonable approach as our main targets are area specific parameters and {\it not} the global parameter $A$.  Obviously, if $A$ is the main target, we would recommend a standard variance component method.
We stress that in general {\it none} of the existing methods for estimating $A$ satisfies any of all the four properties simultaneously.

In Section 2, we propose a new adjusted maximum likelihood estimator of $A$ that satisfies all the four desirable properties listed above. The balanced case has been heavily studied in the literature.  We consider the balanced case in Section 3 and show how our results are related to the ones in the literature.  In Section 4, using a real life data from the U.S. Census Bureau, we demonstrate superior performances of our proposed estimators and MSE estimators over the competing estimators.  A Monte Carlo simulation study, described in Section 5, shows that the proposed estimators outperform competing estimators.  All the technical proofs are deferred to the Appendix.

\section{A New Adjusted Maximum Likelihood Estimator of $A$}
The residual maximum likelihood estimator of $A$ is defined as:
$$\hat A_{RE}=\argmax_{A \in [0,\infty)}L_{RE}(A),$$
where $L_{RE}(A)$ is the residual likelihood of $A$ given by
$$L_{RE}(A)=|X'V^{-1}X|^{-\frac{1}{2}}|V|^{-\frac{1}{2}}\exp\left (-\frac{1}{2}y^{\prime}Py\right ),$$
with $P=V^{-1}-V^{-1}X(X^{\prime}V^{-1}X)^{-1}X^{\prime}V^{-1}$.
Note that $\hat A_{RE}$ does not satisfy any of the four desirable properties listed in the introduction.

In an effort to find a likelihood-based estimator of $A$ that satisfies all the four desirable properties, we define the followed adjusted maximum likelihood estimator of $A$:
$$\hat A_i=\argmax_{A \in [0,\infty)}h_i(A)L_{RE}(A),$$
where $h_i(A)$ is a factor to be suitably chosen so that all the four desirable properties are satisfied.

We first find $h_i(A)$ so that the resulting estimator of $A$ results in a nearly unbiased estimator of $B_i$ that also protects EBLUP from overshrinking.  In other words, we  first find the adjustment factor $h_i(A)$ that simultaneously satisfies Properties 1 and 2.  Interestingly, it turns out that such a adjusted maximum likelihood estimator also satisfies Properties 3 and 4.


Using Lemma \ref{Aihat} in Appendix A and Taylor series expansion, we have
\begin{align}
\mbox{Var} (\hat B_i )=\frac{2D_i^2}{(A+D_i)^4\mbox {tr}[V^{-2}]}+o(m^{-1}),\label{varB}
\end{align}
for large $m$.  We restrict ourselves to the class of estimators of $A$ that satisfies (\ref{varB}).

Using Lemma \ref{Aihat} and Taylor series expansion, we have
\begin{eqnarray}
E(\hat B_i)&=&B_i+\left [ \frac{\partial B_i}{\partial A} \frac{\partial \log h_i(A)} {\partial A}+\frac{1}{2}\frac{\partial^2 B_i}{\partial A^2}\right ]\frac{2}{tr[V^{-2}]}+o(m^{-1}).\label{biasB}
\end{eqnarray}

Thus, Property 1 is satisfied if we have
\begin{align}
\frac{\partial B_i}{\partial A}\frac{\partial \log h_i(A)} {\partial A}+\frac{1}{2}\frac{\partial^2 B_i}{\partial A^2}=0. \label{diff1}
\end{align}

%

Now the differential equation (\ref{diff1}) simplifies to:
\begin{align}
\frac{\partial \log h_i(A)} {\partial A}=\frac{1}{A+D_i}.\label{diff3}
\end{align}
Thus, an adjustment factor that satisfies (\ref{diff3}) is given by
\begin{eqnarray*}
h_{i0} (A)=(A+D_i).
\end{eqnarray*}
This adjustment factor is indeed the unique solution to (\ref{diff1}) up to the order $O(m^{-1})$.  Let $\hat A_{i0}$ be the adjusted maximum likelihood estimator of $A$ for the choice $h_i(A)=h_{i0}(A).$  We note that $\hat A_{i0}$ is not strictly positive.  To achieve strict positivity, we propose our final estimator of $A$ as:
$$\hat {A}_{i;MG}=\argmax_{A \in [0,\infty)}\tilde h_i(A)L_{RE}(A),$$
where $\tilde h_i(A)=h_{+}(A)h_{i0}(A)$ with the additional adjustment $h_{+}(A)$ satisfying regularity conditions R4 and R6-R7.

Our proposed estimator of $B_i$ and EBLUP are given by
$$\hat B_{i;MG}=B_i(\hat A_{i;MG}),\;\;\hat \theta_{i;MG}^{EB}=\hat\theta_i^{BLUP}(\hat A_{i;MG}),$$
respectively.

Unlike the common practice, we avoid bias correction in obtaining both Taylor series and parametric bootstrap MSE estimators of our proposed EBLUP.  Interestingly, our approach ensures the important dual property of MSE estimator --- second-order unbiasedness and strict positivity.  This kind of MSE estimators is the first of its kind in the small area estimation literature.

We obtain our Taylor series estimator of MSE of EBLUP by simply plugging in the proposed estimator $\hat A_{i;MG}$ for $A$ in the second-order MSE approximation $M_{approx}(A)$ and is given by:
\begin{align}
\hat M_{i;MG}\equiv M_{i;approx}(\hat A_{i;MG})=g_{1i}(\hat A_{i;MG})+g_{2i}(\hat A_{i;MG})+g_{3i}(\hat A_{i;MG}),\label{mseTaylor}
\end{align}

Our proposed parametric bootstrap MSE estimator retains the simplicity of bootstrap originally intended in Efron (1979).  It is given by
\begin{align}
\hat{M}_{i;MG}^{boot}\equiv E_{*}[\hat{\theta}_i(\hat A_{i;MG}^*,y^{*})-{\theta}_i^{*}]^2,\label{mseboot}
\end{align}
where ${\theta}_i^{*}=x_i^{\prime}\hat{\beta}(\hat A_{1;MG},\cdots,\hat A_{m;MG})+v^{*}_i$ with
$v^{*}_i\sim N(0, \hat A_{i;MG})$.
Note that the new bootstrap MSE estimator does not require any bias correction.

The following theorem states that our proposed adjusted maximum likelihood estimator of $A$ satisfies all the four desirable properties.
\begin{thm}
\label{multiplegoal}
Under the regularity conditions $R1-R7$, we have, for large $m$,
\begin{description}

\item (i)$Bias( \hat B_{i;MG})=o(1);\;Var( \hat B_{i;MG} )=\frac{2D_i^2}{(A+D_i)^4 \mbox{tr} [V^{-2}]}+o(m^{-1});$

\item (ii) $0<\mbox{inf}_{m\ge 1}\hat B_{i;MG}\le \mbox{sup}_{m\ge 1}\hat B_{i;MG}<1$, for $m>p+2$;

\item (iii)$E( \hat M_{i;MG} ) - M_{i}( \hat {\theta}_{i;MG}^{EB} )=o(m^{-1})$;

\item (iv) $E (\hat{M}_{i;MG}^{boot} ) -  M_{i} (\hat {\theta}_{i;MG}^{EB} ) =o(m^{-1}).$

\end{description}
\end{thm}
For proof of Theorem \ref{multiplegoal}, see Appendix B.



\section{The balanced case:  $D_i=D,\;i=1,\cdots,m$}
In this section, we show how the proposed adjusted maximum likelihood estimator of $A$ is related to the problem of simultaneous estimation of several independent normal means, a topic for intense research activities, especially in the 60's, 70's and 80's, since the introduction of the celebrated James-Stein estimator (James and Stein 1961).

Let $y_i|\theta_i\stackrel{ind}\sim N(\theta_i,1),\;i=1,\cdots,m$.  James and Stein (1961) showed that for $m\ge 3$ the maximum likelihood (also unbiased) estimator of $\theta_i$ is inadmissible under the sum of squared error loss function $L(\hat\theta,\theta)=\sum_{j=1}^m(\hat\theta_j-\theta_j)^2$ and is dominated by the James-Stein estimator:  $\hat\theta_i^{JS}=(1-\hat B_{JS})y_i$, where $\hat B_{JS}={(m-2)}/{\sum_{j=1}^my_j^2}.$  That is,
\begin{eqnarray}
E\left [\sum_{j=1}^m(\hat\theta_j^{JS}-\theta_j)^2|\theta\right ]\le E\left [\sum_{j=1}^m(y_j-\theta_j)^2|\theta\right ],\;\forall \theta\in R^m,\label{eq:js}
\end{eqnarray}
where $R^m$ is the $m$-dimensional Euclidean space, with strict inequality holding for at least one point $\theta$.  The dominance result, however, does not hold for individual components.

Efron and Morris (1973) offered an empirical Bayesian justification of the James-Stein estimator under the prior $\theta_i\stackrel{iid}\sim N(0,A),\;i=1,\cdots,m.$  Their model is indeed a special case of two level normal hierarchical model with $D_i=1,\;x_i^{\prime}\beta=0,\; i=1,\cdots,m,$ and thus the James-Stein estimator of $\theta_i$ can be also viewed as an EBLUP.

Morris (1983) discussed an empirical Bayesian estimation of $\theta_i$ for a Bayesian model that is equivalent to the balanced case of NHM, that is, when $D_i=D$ implying $B_i=B,\;i=1,\cdots,m.$  In this case, he noted that $\hat B_U={(m-p-2)D}/{S}$   is an  exact unbiased estimator of $B$, using the fact that, under NHM, $S=\sum_{j=1}^m(y_j-x_j^{\prime}\hat\beta_{ols})^2\sim (D+A)\chi^2_{m-p},$ where $\hat\beta_{ols}$ is the ordinary least square estimator of $\beta$.  We can write  $\hat B_U\equiv  B(\hat A_{Morris})= {D}/{(D+\hat A_{Morris})}$, where $\hat A_{Morris}={S}/{(m-p-2)}-D$.  One can alternatively estimate $B$ by a simple plug-in estimator: $\hat B_{plug}\equiv B(\hat A_U)={D}/{(D+\hat A_{U})}$, where
$\hat A_U={S}/{(m-p)}-D$ is an  unbiased estimator of $A$.
Note that for $m>p+4$
\begin{eqnarray*}
&&E(\hat B_U-B)=0,\;\;E(\hat B_{plug}-B)=\frac{2}{m-p-2}B=O(m^{-1}),\\
&&V(\hat B_U)=\left (\frac {m-p-2}{m-p}\right )^2 V(\hat B_{plug})\le V(\hat B_{plug}).
\end{eqnarray*}
Thus, $\hat B_U$ is better than $\hat B_{plug}$ both in terms of bias and variance properties.  We can write $\hat B_U=\hat B_{plug}{(m-p-2)}/{(m-p)}.$
As pointed out by Morris (1983), the factor $(m-p-2)/(m-p)$ helps correct for the curvature dependence of $B$ on $A$.

Consider the following empirical Bayes estimator (same as EBLUP) of $\theta_i$:
\begin{eqnarray}
\hat\theta_i^{EB}(\hat A_{Morris})=(1-\hat B_U)y_i+\hat B_U x_i^{\prime}\hat\beta_{ols}.\label{eb}
\end{eqnarray}
In this case, exact MSE and exact unbiased estimator of MSE can be obtained.  Componentwise, for $m\ge p+3$,  we have
$$E[(\hat\theta_i^{EB}(\hat A_{Morris})-\theta_i)^2]\le D.$$

Thus, $\hat\theta_i^{EB}(\hat A_{Morris})$ dominates $y_i$ in terms of unconditional MSE for $m\ge p+3$.
Such a componentwise dominance property, however, does not hold for conditional MSE (conditional on $\theta$); see Morris (1983) for details.

Since $B<1$, using Stein's argument,   Morris (1983) suggested the following estimator of $B$ : $\hat B_{Morris}={D}/{(D+\hat A_{Morris}^{+})},$ where
$\hat A_{Morris}^{+}={S}/{(m-p-2)}-D$ if $S> (m-p-2)D$ and $\hat A_{Morris}^{+}={2D}/{(m-p-2)}$ otherwise.
This improves the estimation of both $B$ and $\theta_i$.  It is straightforward to show that in this special case $\hat A_{Morris}^{+}$  satisfies all the four properties.  Moreover, under the regularity condition R6-R8 and $m>p+2$, $\hat A_{MG}$, our proposed  estimator of $A$, is unique (see Appendix C for a proof) and is equivalent to $\hat A_{Morris}^{+}$ in the higher-order asymptotic sense, that is, $E(\hat A_{MG}-\hat A_{Morris}^{+})=o(m^{-1})$.

Let $\hat \theta_i^{EB}=\hat \theta_i^{EB}(\hat A)$ denote an EBLUP of $\theta_i$, where $\hat A$ could be $\hat A_{MG},\; \hat A_{Morris}^{+}$ or the REML $\hat A_{RE}=\mbox{max}(0,\hat A_U).$  We can write $M_{i;approx}(A)=g_1(A)+g_2(A)+g_3(A)$ as the second-order approximation to  $M_i(\hat \theta_i^{EB})=MSE(\hat \theta_i^{EB})$ for any of the three choices of the estimator of $A$.  The traditional second-order unbiased MSE estimator is obtained by correcting bias of $M_{i;approx}(\hat A_{RE})$, up to the order $O(m^{-1})$.  It is given by  $\hat {M}_{i,RE}=g_1(\hat A_{RE})+g_2(\hat A_{RE})+2g_3(\hat A_{RE})$; see Prasad and Rao (1990), Datta and Lahiri (2000), Das et al. (2004).  In this paper, we suggest an alternative second-order unbiased MSE estimator without bias-correction, that is,  $\hat M_{i;MG}=g_1(\hat A_{MG})+g_2(\hat A_{MG})+g_3(\hat A_{MG})$.

We can show that
\begin{eqnarray*}
V(\hat {M}_{i,RE})&=&a_m+o(m^{-1}),\\
V[\hat M_{i;MG}]&=&b_m+o(m^{-1}),
\end{eqnarray*}
where
\begin{eqnarray*}
a_m&=&\left [\frac{(m-4-mq_i)(m-p)}{m(m-p-2)}\right ]^2\frac{2D^2B^2}{m-p-4},\\
b_m&=&\left (\frac{m-2-mq_i}{m}\right )^2\frac{2D^2B^2}{m-p-4}, \ q_i=x_i^{\prime}(X^{\prime}X)^{-1}x_i.
\end{eqnarray*}
It is straightforward to check that for $m>p+4$ and $p\ge 3$, $b_m\le a_m$.  Thus, in the higher-order asymptotic sense, $\hat M_{i;MG}$ is a better second-order unbiased estimator of
$M_i(\hat \theta_i^{EB})$ than $\hat {M}_{i,RE}.$

\section{A Connection to the Bayesian Approach}
In this section, we suggest a Bayesian method that is close to our proposed EBLUP in certain higher-order asymptotic sense.  To this end, we seek a prior on the hyperparameters $(\beta,A)$ that  satisfies all the  following properties simultaneously:
\begin{description}
\item[(i)] $E[B_i|Y=y]=\hat{B}_{i,MG}+o_p(m^{-1})$;
\item[(ii)] $V[B_i|Y=y]=Var( \hat B_{i;MG} )+o_p(m^{-1})$;
\item[(iii)] $E[\theta_i|Y=y]=\hat{\theta}_{i,MG}+o_p(m^{-1})$;
\item [(iv)] $V[\theta_i|Y=y]=\hat{M}_{i,MG}+o_p(m^{-1})$;
\item [(v)] $V[\theta_i|Y=y]=\hat{M}_{i;MG}^{boot}+o_p(m^{-1})$.
\end{description}
First assume the following prior for $(\beta,A)$:  $p (\beta, A)\propto \pi (A),\;\beta\in R^p,\;A>0$.  We first find a prior $\pi (A)$ satisfying property (i).  To this end, following Datta et al. (2005), we first introduce the following notations:
\begin{align*}
\hat b_1&=\frac{\partial B_i}{\partial A}\Big |_{\hat{A}_{RE}},\;\;\; \hat b_2=\frac{\partial^2 B_i}{\partial A^2}\Big |_{\hat{A}_{RE}}, \ \
\hat \rho_1 =\frac{\partial \log \pi(A)}{\partial A}\Big |_{\hat{A}_{RE}},\\
\hat h_2&=-\frac{1}{m}\frac{\partial^2 l_{RE}}{\partial A^2}\Big |_{\hat{A}_{RE}}=\frac{tr[V^{-2}]}{2m}+o_p(m^{-1}),\\
\hat h_3&=-\frac{1}{m}\frac{\partial^3 l_{RE}}{\partial A^3}\Big |_{\hat{A}_{RE}}=-\frac{2tr[V^{-3}]}{m}+o_p(m^{-1}),
\end{align*}
where $\hat A_{RE}$ is the residual maximum likelihood estimator of $A$, $l_{RE}$ is the logarithm of residual likelihood, and $V=\mbox{diag}(A+D_1,\cdots,A+D_m).$

We have
\begin{align}
\hat{B}_i(\hat{A}_{i,MG})-\hat{B}_i(\hat{A}_{RE})&=(\hat{A}_{i,MG}-\hat{A}_{RE})\hat b_1+\frac{1}{2}(\hat{A}_{i,MG}-\hat{A}_{RE})^2 \hat b_2+O_p(m^{-2})\notag\\
&=(\hat{A}_{i,MG}-A)\hat b_1+o_p(m^{-1})=-\frac{2D_i}{tr[V^{-2}](A+D_i)^3}+o_p(m^{-1});\label{bias.B}
\end{align}
\begin{align}
E[B_i|Y=y]=\hat{B}_i(\hat{A}_{RE})+\frac{1}{2m\hat{h}_2}\left(\hat b_2-\frac{\hat{h}_3}{\hat{h}_2}\hat b_1 \right)+\frac{\hat{b}_1}{m\hat{h}_2}\hat \rho_1+O_p(m^{-2}).\label{B.HB}
\end{align}

It is interesting to note that (\ref{B.HB}) is given by (21) in Datta et al. (2005) with $b(A)=B_i(A)$.
Hence, we seek $\rho_1$ satisfying the following differential equation:
\begin{align}
\frac{1}{2m{h}_2}\left( b_2-\frac{{h}_3}{{h}_2} b_1 \right)+\frac{{b}_1}{m{h}_2} \rho_1=-\frac{2D_i}{tr[V^{-2}](A+D_i)^3}.\label{diff1}
\end{align}

The equation (\ref{diff1}) can be written as follows (up to $O_p(m^{-1}))$;
\begin{align}
\rho_1=\frac{\partial \log \pi(A)}{\partial A}&=-\frac{m  h_2}{ b_1}\frac{2D}{tr[V^{-2}](A+D_i)^3}-\frac{1}{2}\left[\frac{ b_2}{ b_1}-\frac{ h_3}{ h_2} \right]\notag\\
&=\frac{2}{A+D_i}-\frac{2tr[V^{-3}]}{tr[V^{-2}]}.\label{difffinal}
\end{align}
A solution to  differential equation (\ref{difffinal}) is given by;
\begin{eqnarray}
\pi(A)\propto (A+D_i)^2tr[V^{-2}]. \label{MG.pri}
\end{eqnarray}
It is straightforward to check that the prior  (\ref{MG.pri}) satisfies rest of the conditions (ii)-(v).  Interestingly, this prior is same as the prior suggested by Datta et al. (2005).  For the balanced case, the prior reduces to the Stein's harmonic prior.

\section{SAIPE data analysis}
For purposes of evaluation, we consider the problem of estimating the percentages of school-age (aged 5-17) children in poverty for the fifty states and the District of Columbia using the same data set considered by Bell (1999). We choose two years (1992 and 1993) of state level data from the U.S. Census Bureau's Small Area Income and Poverty Estimates (SAIPE) program. In 1992, the REML estimate of A is zero while in year 1993 it is positive.  Thus, these years would provide two different scenarios for evaluating estimation methods.

We assume the standard SAIPE state level model in which survey-weighted estimates of the percentages of 5-17-year-old (related) children in poverty  follow the Fay-Herriot model (\ref{FH}). The survey-weighted percentages are obtained using the Current Population Survey (CPS) data with their sampling variances $D_i$ estimated by a Generalized Variance Function (GVF) method, following Otto and Bell (1995).  However, as in any data analysis that use the Fay-Herriot model, we assume the sampling variances to be known throughout the estimation procedure. We use the same state level auxiliary variables $x$ (a vector of length 5, i.e., $p = 5$), obtained from the Internal Revenue Service (IRS) data, food stamp data and census residual data that the SAIPE program used for the problem.

\begin{figure}[!h]
\begin{center}
\includegraphics[width=15cm,bb=0 0 1567 664]{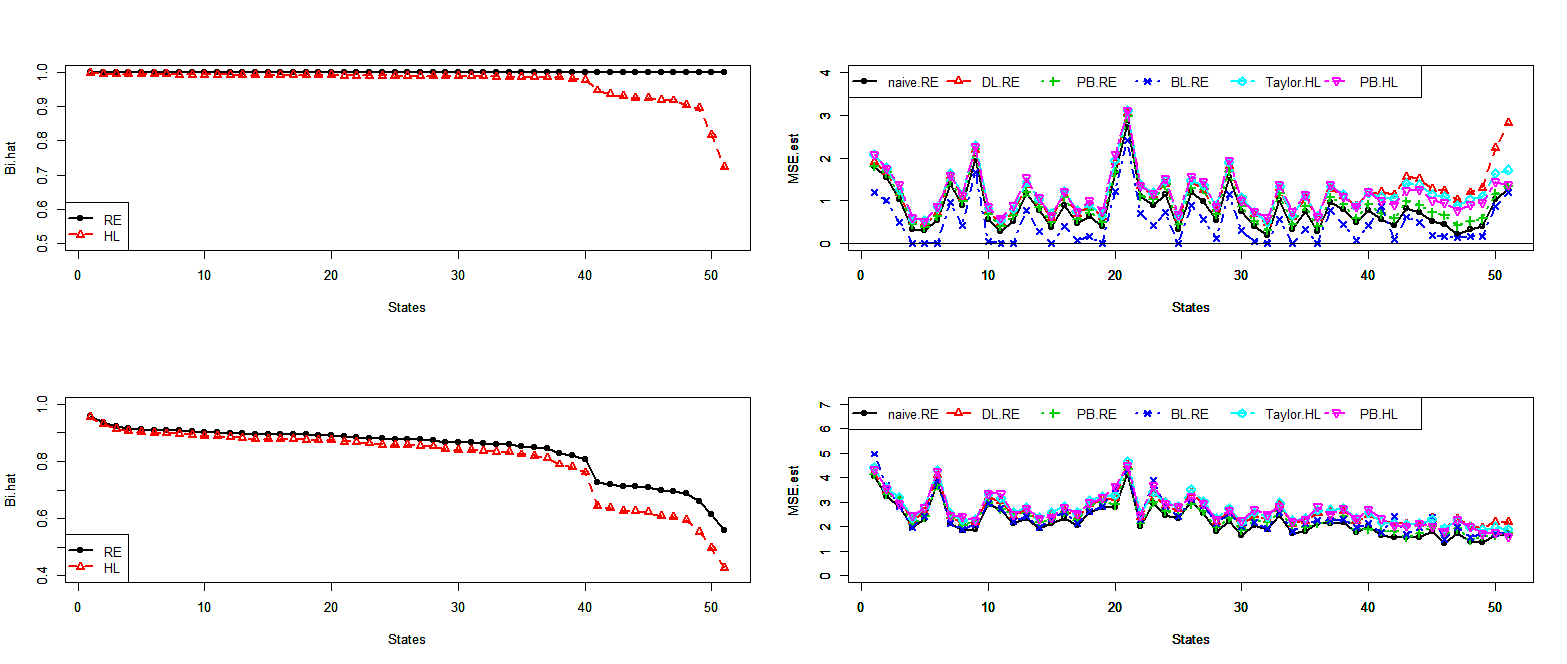}
\caption{Estimates of ${B_i}$ and MSE using all SAIPE data for 1992 (above) and 1993(bottom) year}
\label{saipe.DA.fig}
\end{center}
\end{figure}

Table \ref{B.hat.saipe} displays REML and our proposed estimates (HL) of the shrinkage parameters $B_i$ for Washington DC (DC), Hawaii (HI) and California (CA) for the year 1992 and DC, Oregon (OR) and CA for the year 1993.  They have the largest, median and smallest sampling variances $D_i$ among all the states and DC, respectively.  For 1992, REML estimate of $A$ is zero yielding a $B_i$ estimate of 1 for all the states and DC.  This overshrinkage problem reduces EBLUPs for all the states to regression synthetic estimates.  Thus, even for states  with reliable direct estimates  (e.g., CA), there is no contribution of direct estimates in the EBLUP formula.  Our proposed estimates of shrinkage parameters offer a sensible solution.  For DC, our shrinkage estimate is very close to 1 (giving nearly zero weight to the survey-weighted direct estimate in the EBLUP formula), but for California survey estimate gets considerable weight (about $28\%$).  In 1993, we do not have overshrinkage problem for REML estimates of the shrinkage factors, but our proposed estimates of $B_i$ always gives more weights to the survey-weighted direct estimates than the corresponding REML estimates. Both REML and proposed estimates of $B_i$ for all the states and DC are displayed in the left panel of Figure \ref{saipe.DA.fig}.  Overall, our proposed estimates of $B_i$ are more conservative than REML.

Table \ref{M.hat.saipe} displays different MSE estimates of EBLUPs  for the selected three states for both years.
The right panel of Figure \ref{saipe.DA.fig} displays different MSE estimates for all the states in both years.  For this study, we included the following MSE estimators of EBLUP:
\begin{description}
\item (a) Naive MSE estimator (naive.RE) given by  $g_{1i}(\hat A_{RE})+g_{2i}(\hat A_{RE})$, where $\hat A_{RE}$ denotes the REML estimator of $A$.  This MSE estimator neither incorporates the extra uncertainty due to the estimation of $A$ nor adjusts bias of the estimator $g_{1i}(\hat A_{RE})$ and is not second-order unbiased;

\item (b) Single parametric bootstrap MSE estimator (PB.RE) that is obtained from (\ref{mseboot}) when REML estimator  of $A$ is used in the EBLUP formula and is not a second-order unbiased.

\item (c) Two second-order unbiased MSE estimators based on Taylor-series:
    \begin{description}
    \item (i) DL.RE: $g_{1i}(\hat{A}_{RE})+ g_{2i}(\hat{A}_{RE})+2 g_{3i}(\hat{A}_{RE})$; see  Datta and Lahiri (2000).
    \item (ii) Taylor.HL: the proposed Taylor series MSE estimator given by (\ref{mseTaylor}).
    \end{description}

\item (d) Two second-order unbiased single parametric bootstrap MSE estimators:
     \begin{description}
     \item (i)  BL.RE:   $\;2\{g_{1i}(\hat{A}_{RE})+ g_{2i}(\hat{A}_{RE})\}-E_*[g_{1i}(\hat{A}_{RE}^*)+ g_{2i}(\hat{A}_{RE}^*)] \\
+E_{*}[\{\hat{\theta}_i^{*}(y_i,\hat{A}_{RE}^*,\hat{\beta}(\hat{A}_{RE}^*,y_i))-\tilde{\theta}_i^{*}(y_i,\hat{A}_{RE},\hat{\beta}(\hat{A}_{RE},y_i))\}^2]$;
  see Butar and Lahiri (2003).
     \item (ii) PB.HL: our proposed single parametric bootstrap MSE estimator given by (\ref{mseboot}).
     \end{description}
\end{description}

For this application, there is no appreciable difference between the naive MSE estimates and MSE estimates that attempt to capture additional variability due to the estimation of $A$.  In most of the cases, naive MSE estimates are slightly lower than both the first-order and second-order MSE estimates.   The first-order unbiased MSE estimates (PB.RE) are generally slightly smaller than the second-order unbiased  MSE estimates.  The PB.BL MSE estimates can take negative values because of the adjustment needed to make it second-order unbiased.  Except for large states (e.g., CA),  MSE estimates for EBLUPs are considerably lower than the corresponding sampling variances $D_i$ indicating possible improvements by EBLUPs over the direct estimates.

For the year 1992, REML estimate of $A$ is zero.  This is probably causing unusual behavior for DL.RE or BL.RE MSE estimates.  For example, DL.RE MSE estimate for a large state like CA is more than that for a small state DC (similar behavior can be observed for BL.RE).  For CA, DL.RE MSE estimate is even higher than the corresponding sampling variance of the direct estimate while all the other MSE estimates are showing opposite results.  Overall, our proposed MSE estimates appear reasonable for both years.

\begin{table}[!h]
\small
\caption{Estimates of shrinkage factors $B_i$ in 3 areas (minimum, median and max $B_i$ values) in 1992 and 1993 SAIPE data}
~~\\
\label{B.hat.saipe}
\centering
\begin{tabular}{l|r|rr|l|r|rr}
  \hline
\multicolumn{4}{c|}{1992 year}& \multicolumn{4}{|c}{1993 year}\\ \hline
States & $D_i$ & RE & HL & States & $D_i$ & RE & HL \\
  \hline
DC & 31.6940 & 1.0000 & 0.9968 & DC & 38.2260 & 0.9574 & 0.9546 \\
HI & 11.3470 & 1.0000 & 0.9887 & OR & 12.1880 & 0.8775 & 0.8563 \\
CA & 1.8830 & 1.0000 & 0.7227 & CA & 2.1560 & 0.5588 & 0.4284 \\
   \hline
\end{tabular}
\end{table}

\begin{table}[!h]
\small
\caption{Estimates of MSEs in 3 areas (minimum, median and max $B_i$ values) in 1992 and 1993 SAIPE data}
\label{M.hat.saipe}
\centering
\begin{tabular}{l|r|rrrrrr}
  \hline
\multicolumn{8}{c}{1992 data} \\\hline
 States & $D_i$ & naive.RE & DL.RE & PB.RE & BL.RE & Taylor.HL & PB.HL \\
  \hline
 DC & 31.69 & 1.81 & 1.91 & 1.80 & 1.19 & 2.08 & 2.07 \\
  HI & 11.35 & 1.19 & 1.45 & 1.30 & 0.88 & 1.48 & 1.57 \\
  CA & 1.88 & 1.26 & 2.82 & 1.34 & 1.20 & 1.72 & 1.37 \\ \hline
\multicolumn{8}{c}{1993 data} \\\hline
 States & $D_i$ & naive.RE & DL.RE & PB.RE & BL.RE & Taylor.HL & PB.HL \\
  \hline
 DC & 38.23 & 4.07 & 4.23 & 4.14 & 4.97 & 4.41 & 4.33 \\
  OR & 12.19 & 3.02 & 3.39 & 2.91 & 3.13 & 3.52 & 3.21 \\
  CA & 2.16 & 1.64 & 2.19 & 1.74 & 1.72 & 1.87 & 1.60 \\
   \hline
\end{tabular}
\end{table}

\section{Monte Carlo simulation}
In this section, we report results from a Monte Carlo simulation study.  In particular, we evaluate finite sample performances of two different estimators of $A$ --- the commonly used REML $\hat A_{RE}$ and the proposed  estimator $\hat A_{MG}$ --- in estimating the shrinkage parameters $B_i$, small area means $\theta_i$ and MSE of EBLUPs of $\theta_i$.    To understand the effect of small $m$ on different estimation problems, we set $m=15$ and generate $\{(y_i,\theta_i),\;i=1,\cdots,m\}$ using the Fay-Herriot model (\ref{FH}).

We use the 1992 SAIPE data described in the previous section to design our simulation study.  The 15 areas correspond to states with largest sampling variances $D_i$.  In the simulation, we use $x_i$ and $D_i$ for these states from the 1992 SAIPE data and use $A=15.94$, which is the median of $D_i$ for the 15 states.  The weighted least squared estimates of $\beta$ from the real data including all 50 states and DC are treated as true $\beta$ for the simulation.

We define the relative bias (RB) and relative root mean squared error (RRMSE) of an estimator $\hat B_i$ of $B_i$ as:
\begin{align*}
{\rm RB \ of \ } \hat B_i: \ \frac{{\rm E } (\hat B_i -B_i)}{B_i}\times 100;\\
{\rm RRMSE \ of \ }\hat B_i: \frac{ \sqrt{ {\rm MSE}(\hat B_i)}}{B_i}\times 100,
\end{align*}
where ${\rm MSE}(\hat B_i)={\rm E} (\hat B_i-B_i)^2.$  The expectations in the definitions of RB and RRMSE are approximated by Monte Carlo $1,000$ independent samples from the Fay-Herriot model.
The RB and RRMSE of an estimator $\hat M_i$ of $M_i={\rm MSE}(\hat \theta_i)={\rm E}(\hat\theta_i-\theta_i)^2$, where $\hat \theta_i$ is an estimator of $\theta_i$,  are defined similarly.  For the parametric bootstrap method, we use $1,000$ bootstrap samples.

\begin{figure}[!h]
\begin{center}
\includegraphics[width=15cm,bb=0 0 1561 664]{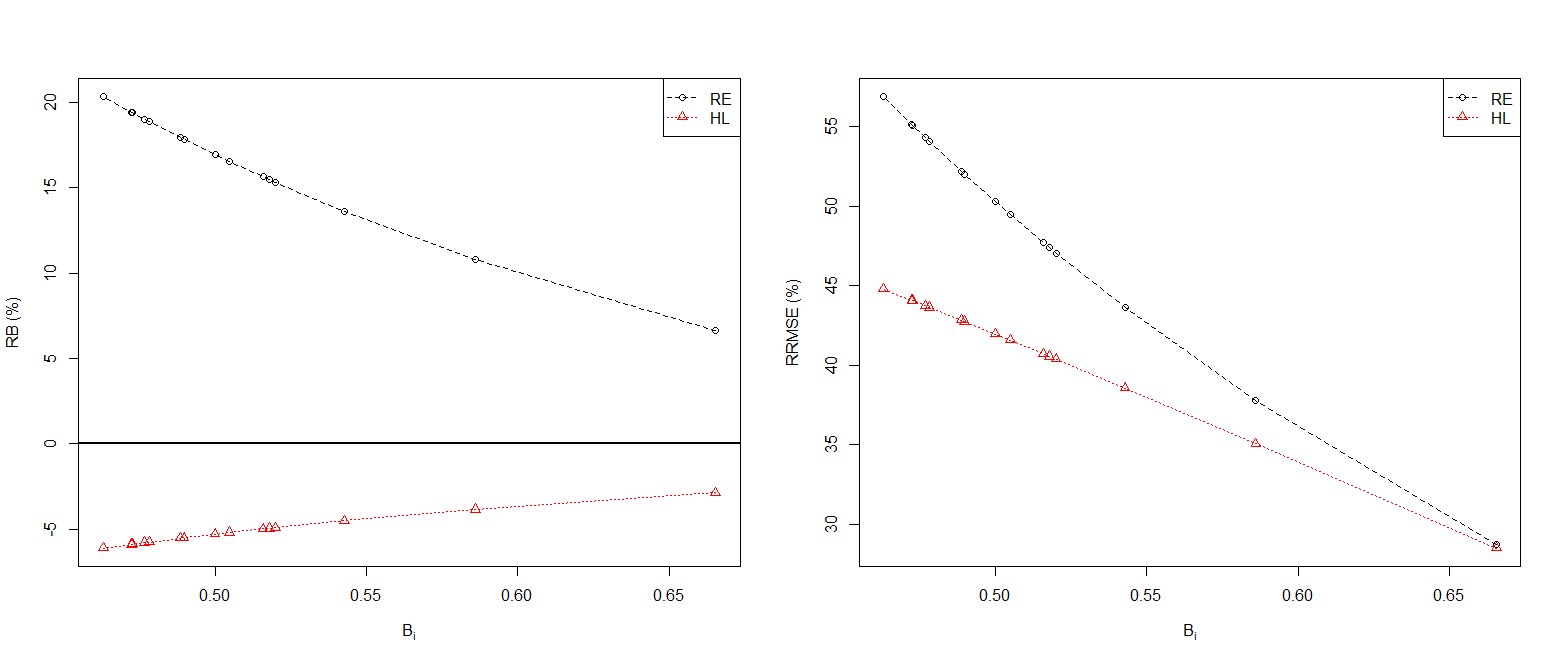}
\caption{RB and RRMSE of $\hat{B}_{i}$}
\label{BhatRB}
\end{center}
\end{figure}

Table \ref{RB.RRMSE.Bhat.tab}  displays simulated  RBs and  RRMSEs of two estimators of $B_i$ for three selected states:  DC, North Dakota (ND), Oklahoma (OK) corresponding to maximum, median and minimum values of $D_i$.  These three states correspond to the maximum (0.67), median (0.50) and minimum values (0.46) of $B_i$'s among the 15 states.  The two estimators of $B_i$ are simple plug-in estimators -- one obtained from REML $\hat A_{RE}$ (denoted by RE) and the other from the proposed  estimator $\hat A_{MG}$ (denoted by HL).  For these states, RE consistently overestimates $B_i$ while HL underestimates.  The absolute values of the RB for HL are always smaller than those of RE.  Moreover, variation of RBs for different $B_i$ is much lower than that of RE.  In terms of RRMSE, HL outperforms RE, especially for small values of $B_i$.
Figure \ref{BhatRB} displays the RB and RRMSE behavior for RE and HL for all the 15 selected states demonstrating superiority of HL over RE.

\begin{figure}[!h]
\begin{center}
\includegraphics[width=8cm,bb=0 0 665 664]{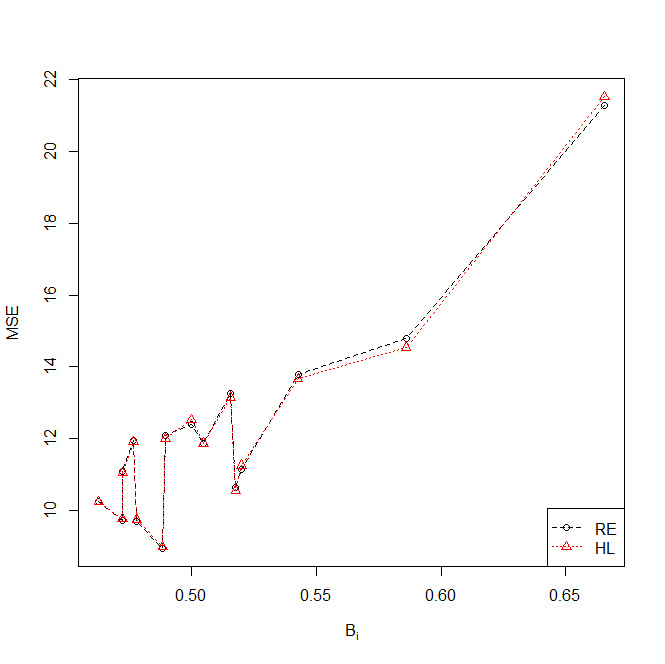}
\caption{MSE of EBLUP with RE and HL }
\label{MSE.EB}
\end{center}
\end{figure}

Figure \ref{MSE.EB} displays the simulated MSEs of two EBLUPs of $\theta_i$ for each of the 15 states, where two EBLUPs are obtained using the REML $\hat A_{RE}$ (RE in the figure) and  estimator $\hat A_{MG}$ (HL in the figure).  There is hardly any difference between the simulated MSEs of the two EBLUPs supporting the theory that these two MSEs are identical up to the order $O(m^{-1}).$

\begin{figure}[!h]
\begin{center}
\includegraphics[width=15cm,bb=0 0 1567 664]{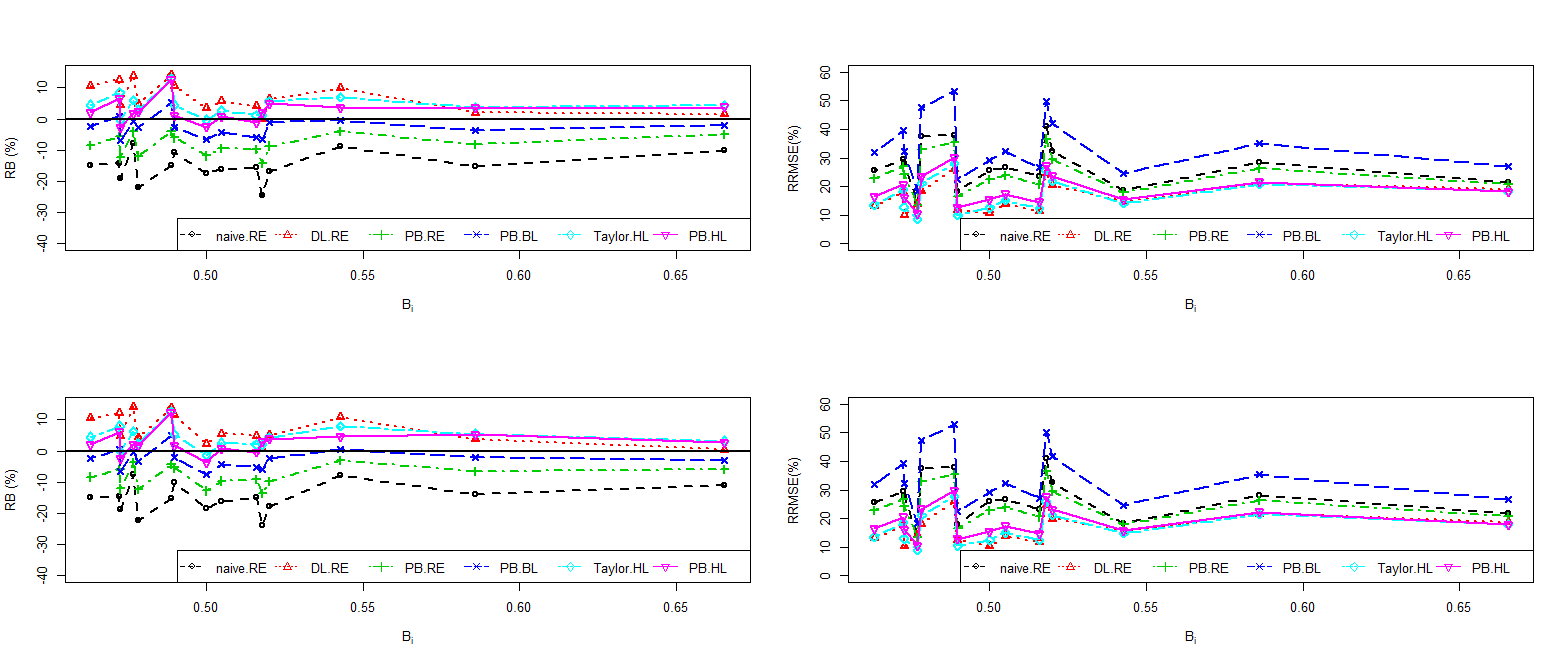}
\label{Mhat.fig}
\caption{ RB and RRMSE of MSE estimators for MSE of EBLUP using REML(above) and HL(bottom)}

\end{center}
\end{figure}

Table \ref{Mhat.RE.tab} reports simulated RBs and RRMSEs of different MSE estimators of EBLUP that uses REML estimator of $A$.  As mentioned earlier, all MSE estimators  except naive.RE and PB.RE are second-order unbiased.  The naive estimator naive.RE consistently underestimates.  All the other MSE estimators improve on naive.RE.  The parametric bootstrap estimator PB.RE that uses REML and does not use bias correction continues to underestimate.  The second-order unbiased parametric bootstrap MSE estimator PB.BL that uses bias correction also underestimates although the amount of underestimation is generally smaller than that of PB.RE.  The proposed second-order unbiased MSE estimators --- Taylor.HL and PB.HL --- are quite competitive to the second-order unbiased Taylor series  MSE estimator, DL.RE, which overestimates for the state with smallest $D_i$.  Our single parametric bootstrap second-order unbiased MSE estimator (PB.HL) that does not involve any bias correction is remarkably better than single parametric bootstrap MSE PB.RE (without bias correction) and even second-order unbiased parametric bootstrap MSE estimator PB.BL (with bias correction).  All MSE estimators except PB.BL have lower RRMSE than naive.RE.  It is interesting to note that the second-order unbiased  PB.BL has more RRMSE than naive.RE for all the three states.  This is probably due to the poor performance of REML of $A$ that PB.BL uses.  The REML of $A$ produces zero estimates $12.4\%$ of the times although true $A$ is 15.94.  The performances of DL.RE, Taylor.HL and PB.HL are similar and all are better than PB.RE.  The performances of the MSE estimators of EBLUP using the proposed estimator of $A$ is similar to the results of Table \ref{Mhat.RE.tab}; see  Table \ref{Mhat.HL.tab}.
The RB and RRMSE behavior of all the MSE estimators for all the 15 states are given in Figure \ref{Mhat.fig}.


\begin{table}[!h]
\small
\caption{RB and RRMSE of $\hat{B}_i$ in 3 areas (min, median and max $B_i$)}
\label{RB.RRMSE.Bhat.tab}
\centering
\begin{tabular}{r|r|rr|rr}
  \hline
&&\multicolumn{2}{c|}{RB}&\multicolumn{2}{c}{RRMSE} \\\cline{3-6}
States&$B_i$ & RE & HL & RE & HL \\
\hline
DC& 0.67 & 6.64 & -2.86 & 28.70 & 28.49 \\
ND& 0.50 & 16.95 & -5.28 & 50.29 & 41.96 \\
OK& 0.46 & 20.31 & -6.09 & 56.90 & 44.79 \\
   \hline
\end{tabular}
\end{table}

\begin{table}[!h]
\small
\caption{RB and RRMSE of $\hat{M}_i$ for MSE of EBLUP with REML in 3 areas (min, median and max $B_i$)}
\label{Mhat.RE.tab}
\centering
\begin{tabular}{r|r|rrrrrr}
  \hline
\multicolumn{8}{c}{RB} \\\hline
States&$B_i$ & naive.RE & DL.RE & PB.RE & PB.BL & Taylor.HL & PB.HL \\
  \hline
DC&0.67 & -10.10 & 1.52 & -4.90 & -2.01 & 4.31 & 3.83 \\
ND&0.50 & -17.50 & 3.39 & -11.81 & -6.57 & -0.35 & -2.63 \\
OK&0.46 & -14.94 & 10.48 & -8.41 & -2.51 & 4.43 & 1.96 \\ \hline
\multicolumn{8}{c}{RRMSE} \\\hline
States&$B_i$ & naive.RE & DL.RE & PB.RE & PB.BL & Taylor.HL & PB.HL \\\hline
DC&0.67 & 21.33 & 19.07 & 20.60 & 26.88 & 18.33 & 18.30 \\
ND&0.50 & 25.51 & 10.64 & 22.54 & 29.28 & 12.57 & 15.48 \\
OK&0.46 & 25.68 & 13.07 & 22.91 & 31.91 & 13.52 & 16.47 \\
   \hline
\end{tabular}
\end{table}

\begin{table}[!h]
\small
\caption{RB and RRMSE of $\hat{M}_i$ for MSE of EBLUP with HL in 3 areas (min, median and max $B_i$)}
\label{Mhat.HL.tab}
\centering
\begin{tabular}{r|r|rrrrrr}
  \hline
\multicolumn{8}{c}{RB} \\\hline
States&$B_i$ & naive.RE & DL.RE & PB.RE & PB.BL & Taylor.HL & PB.HL \\
  \hline
DC&0.67 & -11.09 & 0.40 & -5.95 & -3.09 & 3.16 & 2.68 \\
ND& 0.50 & -18.39 & 2.27 & -12.76 & -7.57 & -1.43 & -3.68 \\
OK& 0.46 & -14.91 & 10.51 & -8.38 & -2.48 & 4.46 & 1.99 \\ \hline
\multicolumn{8}{c}{RRMSE} \\\hline
States&$B_i$ & naive.RE & DL.RE & PB.RE & PB.BL & Taylor.HL & PB.HL \\\hline
DC&  0.67 & 21.64 & 18.81 & 20.66 & 26.68 & 17.90 & 17.90 \\
ND&  0.50 & 25.98 & 10.23 & 22.88 & 29.22 & 12.51 & 15.53 \\
OK&  0.46 & 25.67 & 13.10 & 22.91 & 31.92 & 13.54 & 16.48 \\
   \hline
\end{tabular}
\end{table}

\section{Concluding Remarks}  In this paper, we have solved a set of important problems for the well-known Fay-Herriot small area model through a suitably devised adjusted maximum likelihood estimator of the model variance parameter.  We have demonstrated the superiority of our methods over the existing methods analytically and through data analysis and Monte Carlo simulations.  

Can we extend our results to a general linear mixed model? To answer this question, let us consider the following nested error regression model (NERM) considered by Battese et al. (1988): 
\begin{eqnarray}
y_{ij}=\theta_{ij}+e_{ij}=x_{ij}^{\prime}\beta+v_i+e_{ij}, \ (i=1,\ldots,m;\ j=1,\ldots,n_i),\label{NERM}
\end{eqnarray}
where $\{v_1\ldots,v_m\}$ and $\{e_1,\ldots, e_m\}$ are independent with $v_i{\sim}N(0,\sigma_v^2)$ and $e_i{\sim}N(0,\sigma_e^2)$; $x_{ij}$ is a $p$-dimensional vector of known auxiliary variables; $\beta\in R^p$ is a $p$-dimensional vector of unknown regression coefficients; $\psi=(\sigma_v^2, \sigma_e^2)^{\prime}$ is an unknown variance component vector. $n_i$ is the number of observed unit level data in $i$-th area.


The condition for achieving Property 1, we need to solve the following differential equations with shrinkage factor $B_i=\sigma_e^2/(n_i\sigma_v^2+\sigma_e^2)$, under certain regularity conditions: 
\begin{align}
\left[\frac{\partial \log h_{i}(\psi)}{\partial \psi}\right]^{\prime}I_F^{-1}\left[\frac{\partial B_{i}(\psi)}{\partial \psi}\right]=&H(\psi),
\end{align}
where $$\frac{\partial \log h_{i}(\psi)}{\partial \psi}=\left(\frac{\partial \log h_{i}(\psi)}{\partial \sigma_v^2}, \frac{\partial \log h_{i}(\psi)}{\partial \sigma_e^2}\right)^{\prime},$$
$$H(\psi)=-\frac{1}{2}tr\left[\frac{\partial^2 B_{i}(\psi)}{\partial \psi^2}I_F^{-1}\right], \ \
\frac{\partial B_{i}(\psi)}{\partial \psi}=\frac{n_i}{(n_i\sigma_v^2+\sigma_e^2)^2}(-\sigma_e^2, \sigma_v^2)^{\prime},$$ 
$$I_F^{-1}=\frac{2}{a}\left (
\begin{array}{cc}
\sum[(n_i-1)/\sigma_e^4+(n_i\sigma_v^2+\sigma_e^2)^{-2}] & -\sum n_i/(n_i\sigma_v^2+\sigma_e^2)^2\\
-\sum n_i/(n_i\sigma_v^2+\sigma_e^2)^2&\sum n_i^2/(n_i\sigma_v^2+\sigma_e^2)^2\\
\end{array}\right),$$
$$a=[\sum n_i^2/(n_i\sigma_v^2+\sigma_e^2)^2][\sum \{(n_i-1)/\sigma_e^4+(n_i\sigma_v^2+\sigma_e^2)^{-2}\}]-[\sum n_i/(n_i\sigma_v^2+\sigma_e^2)^2]^2.$$

If we use the following adjustment factor for achieving Property 1:
\begin{align}
\frac{\partial \log h_{i}(\psi)}{\partial \psi}=v{k},
\end{align}
with some fixed two dimensional vector ${k}$, 
the solution of $v$ can be obtained as $v=\frac{H(\psi)}{{k}^{\prime}I_F^{-1}\frac{\partial B_{i}(\psi)}{\partial {\psi}}}$ for some $k$. 
This solution thus lead to a suitable adjustment factor satisfying $$\frac{\partial \log h_{i}(\psi)}{\partial \psi}=\frac{H(\psi)}{{k}^{\prime}I_F^{-1}\frac{\partial B_{i}(\psi)}{\partial {\psi}}}{k}.$$ 
Thus, there exists multiple solutions for adjustment factor satisfying Property 1 under NERM. 

To address such a problem, we will search for the most suitable adjustment factor for the general linear mixed model in the future.


\clearpage
\appendix

\def\thesection{\Alph{section}}
\section{Regularity conditions and Lemma 1}

\begin{description}

\item {\it R1:} $\mbox{rank}(X)=p$ is bounded for large $m$;
\item {\it R2:} The elements of $X$ are uniformly bounded, implying $\sup_{j\geq 1}x_j(X^{\prime}X)^{-1}x_j=O(m^{-1})$;
\item {\it R3:} $0<\inf_{j\geq 1}D_j\leq \sup_{j\geq 1}D_j<\infty$, $A\in (0,\infty)$;
\item {\it R4:}  $\log h_i(A)$ is free of $y$ and  four times continuously differentiable with respect to $A$.  Moreover, $\frac{\partial^k \log h_{i}(A)}{\partial A^k}$ is of order $O(1)$, respectively, for large $m$ with $k=0,1,2,3,4$;
\item {\it R5:} $|\hat{A}_{i}|<C_{ad}m^{\lambda}$, where  $C_{ad}$ a generic positive constant and  $\lambda$ is  small positive constant.
\end{description}
In addition to $R4$, the adjustment factor $h_{+}(A)$ satisfy the following regularity conditions:
\begin{description}
\item {\it R6:}  $\log h_{+}(A)$ is free of $y$ and  four times continuously differentiable with respect to $A$.  Moreover, $\frac{\partial^k \log h_{+}(A)} {\partial A^k}$ is of order $o(1)$, for large $m$ with $k=0,1,2,3,4$;
\item {\it R7;}  $h_{+}(A)$ is a strictly positive on $A>0$ satisfying that $h_{+}(A)\Big|_{A=0}=0$ and $h_{+}(A)<C$ on $A>0$ with a generic positive constant $C$;
\item {\it R8:} In balanced case, that is, $D_i=D$ for all $i$, $(A+D)^2\frac{\partial \log h_{+}(A)} {\partial A}$ is a monotonically decreasing function of $A>0$ with $\lim_{A\rightarrow +0}(A+D)^2\frac{\partial \log h_{+}(A)} {\partial A}=\infty$. When we assume that $\frac{\partial \log h_{+}(A)} {\partial A}>0$, then, $\lim_{A\rightarrow \infty}(A+D)^2  \frac{\partial \log h_{+}(A)} {\partial A}=C$ for fixed $m$, where $C$ is a generic positive constant.

\end{description}
The choice of $h_{+}(A)$ is not unique in general.  One can use the choice given in Yoshimori and Lahiri (2014).

We first present the following Lemma that provides properties of $\hat A_i$ of $A$.
The proof of the theorem is immediate from  Theorem 1 of  Yoshimori and Lahiri (2014) and Das et al. (2004).
\begin{lem}
\label{Aihat}
Under the regularity conditions $R1-R5$, we have, for large $m$,
\begin{description}
\item (i) $E(\hat {A}_i-A)=\frac{\partial \log h_i(A)} {\partial A}\frac{2}{{\mbox{tr} [V^{-2}]}} +o(m^{-1});$

\item (ii) $E(\hat {A}_i-A)^2=\frac{2}{\mbox{tr} [V^{-2}]}+o(m^{-1});$

\item (iii) $E [\hat\theta_i^{EB}(\hat A_i)-\theta_i ]^2\equiv M_i[\hat{\theta}_i^{EB}(\hat A_i)]=M_{i;approx}({A})+o(m^{-1})$, where
$M_{i;approx}({A})=g_{1i}(A)+g_{2i}(A)+g_{3i}(A)$ with
$g_{1i}(A)={AD_i}/({A+D_i}),$
$g_{2i}(A)={D_i^2}x_i^{\prime}(X^{\prime}V^{-1}X)^{-1}x_i/{(A+D_i)^2},$
$g_{3i}(A)={2D_i^2}/[(A+D_i)^3 tr\{V^{-2}\}]$.
\end{description}
\end{lem}

\def\thesection{\Alph{section}}
\section{Proofs of  Theorem 1}

\subsection{Proof of part (i)}
First note that the adjustment factor $h_i(A)$ satisfies regularity condition R4.  Then part (i)  follows  from the construction and (\ref{varB}).

\subsection{Proof of part (ii)}
It suffices to show the strictly positivity for $\hat{A}_{i;MG}$. Note that  $h_{+}(A)h_{i0}(A)L_{RE}(A)\Big{|}_{A=0}=0$ and $h_{+}(A)h_{i0}(A)L_{RE}(A)\geq 0$ for $A\ge 0$ using R6-R7.  Thus, we are left to show that $$\lim_{A\rightarrow \infty} h_{+}(A)h_{i0}(A)L_{RE}(A)=0.$$

\indent Let $C$ be a generic constant.  Using regularity conditions and $m\geq 1$, we have
\begin{align*}
h_{+}(A)h_{i0}(A)&< C(A+\sup_{i\geq 1}D_i),\\
L_{RE}(A)&<C(A+\sup_{i\geq 1} D_{i})^{\frac{p}{2}}| X^{\prime}X |^{-\frac{1}{2}}(A+\inf_{i\geq 1}D_{i})^{-\frac{m}{2}},
\end{align*}
which imply
\begin{align*}
0&\leq h_{+}(A)h_{i0}(A)L_{RE}(A)\\
&< C(A+\sup_{i\geq 1}D_{i})^{1+p/2}(A+\inf_{i\geq 1}D_{i})^{-m/2}|X^{\prime}X|^{-1/2}
\approx A^{-\frac{1}{2}(m-p-2)},
\end{align*}
for large $A$. Thus,  $\hat{A}_{i;MG}$ is strictly positive if $m>p+2$.

\subsection{Proof of part (iii)}
Using part (iii) of Lemma \ref{Aihat}, we get
$$M_i(\hat\theta_{i;MG}^{EB})=M_{i;approx}({A})+o(m^{-1}).$$
Note that  using part (i) of Lemma \ref{Aihat} we have:  $E[g_{2i}(\hat A_{i;MG})]=g_{2i}(A)+o(m^{-1}),\;E[g_{3i}(\hat A_{i;MG})]=g_{3i}(A)+o(m^{-1})$.  Since $g_{1i}(A)=(1-B_{i})D_i$, we have $E[g_{1i}(\hat A_{i;MG})]=g_{1i}(A)+o(m^{-1})$,  using part (i).
This proves part (iii).

\subsection{Proof of part (iv)}
Using part (iii), we have
\begin{align*}
\hat{M}_{i;MG}^{boot}=&g_1(\hat {A}_{i;MG})+g_2(\hat {A}_{i;MG})+g_3(\hat {A}_{i;MG})+R,\\
=&{M}_{i}(\hat {A}_{i;MG})+R,
\end{align*}
where $E[|R|]=o(m^{-1})$.
The result now follows from part (iii).

\section{Proof of the uniqueness of $\hat{A}_{MG}$ in balanced case}

In the balanced case, we have
\begin{align*}
\frac{\partial \log L(A)}{\partial A}=\frac{1}{2(A+D)^2}\left[y^{\prime}\{I_m-X(X^{\prime}X)^{-1}X^{\prime}\}y-(m-p)(A+D) \right].
\end{align*}
Thus, $(A+D)^2\frac{\partial \log L(A)}{\partial A}$ is a linear function of $A$.
Therefore, our estimate of $A$ is obtained as a solution of:
\begin{align}
-(m-p-2)(A+D)+2(A+D)^2\frac{\partial \log h_+(A)}{\partial A}+y^{\prime}\{I_m-X(X^{\prime}X)^{-1}X^{\prime}\}y=0.\label{uniq}
\end{align}
Define $K(A)$ as the left hand of (\ref{uniq}).
 For $A>0$, using the regularity condition R6-R8 and $m>p+2$, we show that $\lim_{A\rightarrow +0}K(A)=\infty$, $\lim_{A\rightarrow \infty}K(A)=-\infty$ and $K(A)$ is a strictly monotonically decreasing function of $A$ on $A>0$.
Hence, there exist $A_{+}$ and $A_{-}$ such that $K(A_+)=-\varepsilon$ and $K(A_-)=\varepsilon$ with small $\varepsilon>0$ and $0<A_-<A_+<\infty$.
Thus, using the intermediate value theorem, we conclude that the adjustment term $h_{+}(A)$ leads to a unique estimate of $A$ on $A>0$.

\section*{Acknowledgement}
The first author\rq{}s research was supported by Grant-in-Aid for Research Activity start-up, JSPS Grant Number 26880011.
The second author\rq{}s research was supported in part by the National Science Foundation Grant Number SES-1534413.

\end{document}